\documentclass{amsart}[12pt]

\usepackage{amsmath,amssymb,epsfig,amscd,xy,amsthm}

\newtheorem{theorem}{Theorem}[section]

\newtheorem{lemma}[theorem]{Lemma}

\newtheorem{fact}[theorem]{Fact}
\newtheorem{cor}[theorem]{Corollary}
 
\newtheorem{definition}[theorem]{Definition}

\newtheorem{conjecture}[theorem]{Conjecture}

\theoremstyle{plain}
\numberwithin{equation}{theorem}

\theoremstyle{remark}

\newcommand{\tensor}{\otimes}

\newcommand{\fO}{\mathfrak o}

\newcommand{\Kbar}{\overline{K}}

\DeclareMathOperator{\sep}{sep}

\DeclareMathOperator{\tor}{tor}

\DeclareMathOperator{\CC}{\mathbb{C}_{\infty}}

\newcommand{\bG}{{\mathbb G}}

\newcommand{\bC}{{\mathbb C}}

\newcommand{\bF}{{\mathbb F}}
\newcommand{\Fq}{\bF_q}
\newcommand{\lra}{\longrightarrow}

\newcommand{\bK}{\overline{K}}
\newcommand{\hhat}{{\widehat h}}

\newcommand{\Drin}{{\bf \phi}}

\title{A dynamical version of the Mordell-Lang conjecture for the additive group}
\author{D.~Ghioca}
\email{dghioca@math.mcmaster.ca}
\address{
Dragos Ghioca \\
Department of Mathematics\\
McMaster University \\
1280 Main Street West \\ 
Hamilton, Ontario \\
Canada  L8S 4K1 \\
}

\author{T.~J.~Tucker}
\email{ttucker@math.rochester.edu}
\address{
Thomas Tucker\\
Department of Mathematics\\
Hylan Building\\
University of Rochester\\
Rochester, NY 14627
}

\thanks{The second author was partially supported by National Security
  Agency Grant 06G-067}

\begin{document}

\begin{abstract}
  We prove a dynamical version of the Mordell-Lang conjecture in the context of Drinfeld modules. We use analytic methods similar to the ones employed by Skolem, Chabauty
and Coleman for studying diophantine equations.
\end{abstract}

\maketitle

\section{Introduction}
\label{intro}

Faltings proved the Mordell-Lang conjecture in the following form (see 
\cite{Faltings}).
\begin{theorem}[Faltings]
\label{T:F}
Let $G$ be an abelian variety defined over the field of complex
numbers $\mathbb{C}$. Let $X\subset G$ be a closed subvariety and
$\Gamma\subset G(\mathbb{C})$ a finitely generated subgroup of
$G(\bC)$. Then $X(\mathbb{C})\cap\Gamma$ is a finite union of cosets
of subgroups of $\Gamma$.
\end{theorem}

In particular, Theorem~\ref{T:F} says that an irreducible subvariety
$X$ of an abelian variety $G$ has a Zariski dense intersection with a
finitely generated subgroup of $G(\bC)$ only if $X$ is a translate of
an algebraic subgroup of $G$. We also note that Faltings result was
generalized to semiabelian varieties $G$ by Vojta (see \cite{V1}), and
then to finite rank subgroups $\Gamma$ of $G$ by McQuillan (see
\cite{McQ}).

If we try to formulate the Mordell-Lang conjecture in the context of
algebraic subvarieties contained in a power of the additive group
scheme $\mathbb{G}_a$, the conclusion is either false (in the
characteristic $0$ case, as shown by the curve $y=x^2$ which has an
infinite intersection with the finitely generated subgroup
$\mathbb{Z}\times\mathbb{Z}$, without being itself a translate of an
algebraic subgroup of $\mathbb{G}_a^2$) or it is trivially true (in
the characteristic $p>0$ case, because every finitely generated
subgroup of a power of $\mathbb{G}_a$ is finite). Denis
\cite{Denis-conjectures} formulated a Mordell-Lang conjecture for
powers of $\bG_a$ in characteristic $p$ in the context of Drinfeld
modules. Denis replaced the \emph{finitely generated subgroup} from
the usual Mordell-Lang statement with a \emph{finitely generated
  $\phi$-submodule}, where $\phi$ is a Drinfeld module. He also
strengthened the conclusion of the Mordell-Lang statement by asking
that the \emph{subgroups} whose cosets are contained in the
intersection of the algebraic variety with the finitely generated
$\phi$-submodule be actually \emph{$\phi$-submodules}. The first
author proved several cases of the Denis-Mordell-Lang conjecture in
\cite{IMRN} and \cite{full-ml-drinfeld}.

In the present paper we investigate other cases of the
Denis-Mordell-Lang conjecture through methods different from the ones
employed in \cite{IMRN}.  In particular, we prove the
Denis-Mordell-Lang conjecture in the case where the finitely generated
$\phi$-module is cyclic and the Drinfeld modules are defined over a
field of transcendence degree equal to one (this is our
Theorem~\ref{hypersurface}).  Note that \cite{IMRN} and
\cite{full-ml-drinfeld} treat only the case where the transcendence
degree of the field of definition is greater than one. One of the
methods employed in \cite{IMRN} (and whose outcome was later used in
\cite{full-ml-drinfeld}) was specializations; hence the necessity of
dealing with fields of transcendence degree greater than one.  By
contrast, the techniques used in this paper are more akin to those
used in treating diophantine problems over number fields (see
\cite{Chabauty}, \cite{Coleman}, or \cite[Chapter 4.6]{BS}, for
example), where such specialization arguments are also not available.
So, making a parallel between the classical Mordell-Lang conjecture
and the Denis-Mordell-Lang conjecture, we might say that the papers
\cite{IMRN} and \cite{full-ml-drinfeld} deal with the ``function field
case'', while our present paper deals with the ``number field case''
of the Denis conjecture. Moreover, using specializations (as in \cite{Mordell-Lang} and \cite{IMRN}), our Theorem~\ref{hypersurface} can be extended to Drinfeld modules defined over fields of arbitrary finite transcendence degree.

We also note that recently there has been significant progress on
establishing additional links between classical diophantine results
over number fields and similar statements for Drinfeld modules. The
first author proved in \cite{Mat.Ann} an equidistribution statement
for torsion points of a Drinfeld module, which is similar to the
equidistribution statement established by Szpiro-Ullmo-Zhang
\cite{suz} (which was later extended by Zhang \cite{Zhang} to a full
proof of the famous Bogomolov conjecture). Bosser \cite{Bosser} proved
a lower bound for linear forms in logarithms at an infinite place
associated to a Drinfeld module (similar to the classical result
obtained by Baker \cite{Baker} for usual logarithms, or by David
\cite{David} for elliptic logarithms). Bosser's result was used by the
authors in \cite{findrin} to establish certain equidistribution and
integrality statements for Drinfeld modules. Moreover, Bosser's result
is quite possibly true also for linear forms in logarithms at finite
places for a Drinfeld module. Assuming this last statement, the
authors proved in \cite{siegel-drinfeld} the analog of Siegel's
theorem for finitely generated $\phi$-submodules. We believe that our
present paper provides an additional proof of the fact that the
Drinfeld modules represent the right arithmetic analog in
characteristic $p$ for abelian varieties in characteristic $0$.

The idea behind the proof of our Theorem~\ref{hypersurface} can be
explained quite simply. Assuming that an affine variety
$V\subset\bG_a^g$ has infinitely many points in common with a cyclic
$\phi$-submodule $\Gamma$, we can find then a suitable submodule
$\Gamma_0\subset \Gamma$ whose coset lies in $V$. Indeed, applying the
logarithmic map (associated to a suitable place $v$) to $\Gamma_0$
yields a line in the vector space $\bC_v^g$.  Each polynomial $f$ that
vanishes on $V$, then gives rise to an analytic function $F$ on this
line (by composing with the exponential function).  Because we assumed
there are infinitely many points in $V \cap \Gamma$, the zeros of $F$
must have an accumulation point on this line, which means that $F$
vanishes identically on the line.  This means that there is an entire
translate of $\Gamma_0$ contained in the zero locus of $f$.  The
inspiration for this idea comes from the method employed by Chabauty
in \cite{Chabauty} (and later refined by Coleman in \cite{Coleman}) to
study the intersection of a curve $C$ of genus $g$, embedded in its
Jacobian $J$, with a finitely generated subgroup of $J$ of rank less
than $g$.  Our technique also bears a resemblance to Skolem's method
for treating diophantine equations (see \cite[Chapter 4.6]{BS}).

Alternatively, our results can be interpreted purely from the point of
view of polynomial dynamics, as we describe the intersection of affine
varieties with the iterates of a point in the affine space under
polynomial actions on each coordinate.  In this paper we will treat
the case of affine varieties embedded in $\bG_a^g$, while the
polynomial action (on each coordinate of $\bG_a^g$) will always be
given by Drinfeld modules.  The more general problem of studying
intersections of affine varieties with the iterates of a point in
affine space under polynomial actions over number fields or function
fields appears to be quite difficult.  To our knowledge, very little
about this question has been proven except in the case of
multiplication maps on semiabelian varieties (see \cite{V1} and
\cite{McQ}).  We refer the reader to Section 4 of Zhang's notes
\cite{ZhangLec} for a number of algebraic dynamical conjectures that
would generalize well-known arithmetic theorems for semiabelian
varieties.  Although these notes do not contain a dynamical analog of
the Mordell-Lang conjecture, Zhang has indicated to us that it might
be reasonable to conjecture that if $\psi: Y \lra Y$ is a suitable
morphism of a projective variety $Y$ (one that is ``polarized'', to
use the terminology of \cite{ZhangLec}), then the intersection of the
$\psi$-orbit of a point $\beta$ with a subvariety $V$ must be finite
if $V$ does not contain a positive dimensional preperiodic subvariety.
 
We briefly sketch the plan of our paper. In Section~\ref{notation} we
set the notation, describe the Denis-Mordell-Lang conjecture, and then
state our main result. In Section~\ref{proofs} we prove this main
result (Theorem~\ref{hypersurface}), while in Section~\ref{further
  extensions} we prove a couple of extensions of it
(Theorems~\ref{rational 2} and \ref{full-rational}).

\section{Notation and statement of our main result}
\label{notation}

All subvarieties appearing in this paper are closed.

\subsection{Drinfeld modules}
We begin by defining a Drinfeld module.  Let $p$ be a prime and let
$q$ be a power of $p$. Let $A:=\mathbb{F}_q[t]$, let $K$ be a finite field extension of $\mathbb{F}_q(t)$, and let $\Kbar$ be an
algebraic closure of $K$. Let $K^{\sep}$ be the separable closure of $K$ inside $\Kbar$. We let $\tau$ be the Frobenius on
$\mathbb{F}_q$, and we extend its action on $\Kbar$.  Let $K\{\tau\}$
be the ring of polynomials in $\tau$ with coefficients from $K$ (the
addition is the usual addition, while the multiplication is the
composition of functions).

A Drinfeld module is a morphism $\Drin:A\rightarrow K\{\tau\}$ for
which the coefficient of $\tau^0$ in $\Drin(a)=:\Drin_a$ is $a$ for
every $a\in A$, and there exists $a\in A$ such that $\Drin_a\ne
a\tau^0$. The definition given here represents what Goss \cite{Goss}
calls a Drinfeld module of ``generic characteristic''.

We note that usually, in the definition of a Drinfeld module, $A$ is
the ring of functions defined on a projective nonsingular curve $C$,
regular away from a closed point $\eta\in C$. For our definition of a
Drinfeld module, $C=\mathbb{P}^1_{\mathbb{F}_q}$ and $\eta$ is the
usual point at infinity on $\mathbb{P}^1$. On the other hand, every
ring of regular functions $A$ as above contains $\mathbb{F}_q[t]$ as a
subring, where $t$ is a nonconstant function in $A$.
 
For every field extension $K\subset L$, the Drinfeld module $\Drin$
induces an action on $\mathbb{G}_a(L)$ by $a*x:=\Drin_a(x)$, for each
$a\in A$. We call \emph{$\phi$-submodules} subgroups of
$\mathbb{G}_a(\overline{K})$ which are invariant under the action of
$\phi$. We define the \emph{rank} of a $\phi$-submodule $\Gamma$ be 
$$\dim_{\Fq(t)}\Gamma\tensor_A\Fq(t).$$
If $\phi_1:A\rightarrow
K\{\tau\},\dots,\phi_g:A\rightarrow K\{\tau\}$ are Drinfeld modules,
then $(\phi_1,\dots,\phi_g)$ acts on $\mathbb{G}_a^g$ coordinate-wise
(i.e. $\phi_i$ acts on the $i$-th coordinate). We define as above the
notion of a $(\phi_1,\dots,\phi_g)$-submodule of $\bG_a^g$; same for
its rank.

A point $\alpha$ is \emph{torsion} for the Drinfeld module action if
and only if there exists $Q\in A\setminus\{0\}$ such that
$\Drin_Q(\alpha)=0$. The set of all torsion points is denoted by $\phi_{\tor}$.

\subsection{Valuations}
Let $M_{\mathbb{F}_q(t)}$ be the set of places on $\Fq(t)$.  We denote
by $v_{\infty}$ the place in $M_{\Fq(t)}$ such that
$v_{\infty}(\frac{f}{g})=\deg(g)-\deg(f)$ for every nonzero $f,g\in
A=\Fq[t]$. We let $M_K$ be the set of valuations on $K$. Then $M_K$ is
a set of valuations which satisfies a product formula (see
\cite[Chapter 2]{Serre-Mordell_Weil}). Thus
\begin{itemize}
\item for each nonzero $x\in K$, there are finitely many $v\in M_K$
such that $|x|_v\ne 1$; and \\
\item for each nonzero $x\in K$, we have $\prod_{v\in M_K} |x|_v=1$.
\end{itemize}

\begin{definition}
Each place in $M_K$ which lies over $v_{\infty}$ is called an infinite place. Each place in $M_K$ which does not lie over $v_{\infty}$ is called a finite place.
\end{definition}

By abuse of notation, we let $\infty\in M_K$ denote any place extending the
place $v_{\infty}$.  

For $v\in M_K$ we let $K_v$ be the completion of $K$ with respect to $v$. Let $\mathbb{C}_v$ be the completion of an algebraic closure
of $K_v$.  Then $| \cdot |_v$ extends to a unique absolute value on
all of $\mathbb{C}_v$. We fix
an embedding of $i: \bK  \lra \mathbb{C}_v$.  For $x \in\bK$, we denote
$| i(x) |_v$ simply as $|x|_v$, by abuse of notation.  

\subsection{Logarithms and exponentials associated to a Drinfeld module}
\label{v-adic}
Let $v\in M_K$. According to Proposition $4.6.7$ from \cite{Goss}, there exists an unique formal power series $\exp_{\phi,v}\in\bC_v\{\tau\}$ such that for every $a\in A$, we have 
\begin{equation}
\label{first expos}
\phi_a=\exp_{\phi,v}a\exp_{\phi,v}^{-1}.
\end{equation} 
In addition, the coefficient of the linear term in $\exp_{\phi,v}(X)$ equals $1$. We let $\log_{\phi,v}$ be the formal power series $\exp_{\phi,v}^{-1}$, which is the inverse of $\exp_{\phi,v}$.

If $v=\infty$ is an infinite place, then $\exp_{\phi,\infty}(x)$ is convergent for all $x\in\CC$ (see Theorem $4.6.9$ of \cite{Goss}). There exists a sufficiently small ball $B_{\infty}$ centered at the origin such that $\exp_{\phi,\infty}$ is an isometry on $B_{\infty}$ (see Lemma $3.6$ of \cite{findrin}). Hence, $\log_{\phi,\infty}$ is convergent on $B_{\infty}$. Moreover, the restriction of $\log_{\phi,\infty}$ on $B_{\infty}$ is an analytic isometry (see also Proposition $4.14.2$ of \cite{Goss}).

If $v$ is a finite place, then $\exp_{\phi,v}$ is convergent on a sufficiently small ball $B_v\subset\bC_v$ (this follows identically as the proof of the analyticity of $\exp_{\phi,\infty}$ from Theorem $4.6.9$ of \cite{Goss}). Similarly as in the above paragraph, at the expense of replacing $B_v$ by a smaller ball, we may assume $\exp_{\phi,v}$ is an isometry on $B_v$. Hence, also $\log_{\phi,v}$ is an analytic isometry on $B_v$.

For every place $v\in M_K$, for every $x\in B_v$ and for every polynomial $a\in A$, we have (see \eqref{first expos})
\begin{equation}
\label{logs and expos}
a\log_{\phi,v}(x)=\log_{\phi,v}(\phi_a(x))\text{ and }\exp_{\phi,v}(ax)=\phi_a(\exp_{\phi,v}(x)).
\end{equation}

By abuse of language, $\exp_{\phi,\infty}$ and $\exp_{\phi,v}$ will be called exponentials, while $\log_{\phi,\infty}$ and $\log_{\phi,v}$ will be called logarithms.

\subsection{Integrality and reduction}
\begin{definition}
  A Drinfeld module $\phi$ has good reduction at a place $v$ if for
  each nonzero $a\in A$, all coefficients of $\phi_a$ are $v$-adic
  integers and the leading coefficient of $\phi_a$ is a $v$-adic unit.
  If $\phi$ does not have good reduction at $v$, then we say that
  $\phi$ has bad reduction at $v$.
\end{definition}

It is immediate to see that $\phi$ has good reduction at $v$ if and only if all coefficients of $\phi_t$ are $v$-adic integers, while the leading coefficient of $\phi_t$ is a $v$-adic unit. All infinite places of $K$ are places of bad reduction for $\phi$.

\subsection{The Denis-Mordell-Lang conjecture}
Let $g$ be a positive integer.
\begin{definition}
Let $\phi_1:A\to K\{\tau\},\dots,\phi_g:A\to K\{\tau\}$ be Drinfeld modules. An algebraic $(\phi_1,\dots,\phi_g)$-submodule of $\mathbb{G}_a^g$ is an irreducible algebraic subgroup of $\mathbb{G}_a^g$ invariant under the action of $(\phi_1,\dots,\phi_g)$.
\end{definition}

Denis proposed in Conjecture $2$ of \cite{Denis-conjectures} the
following problem, which we call the \emph{full} Denis-Mordell-Lang
conjecture because it asks for the description of the intersection of
an affine variety with a \emph{finite rank $\phi$-module} (as opposed
to only a finitely generated $\phi$-module).  Recall that a
$\phi$-module $M$ is said to be a finite rank $\phi$-module if
it contains a finitely generated $\phi$-submodule such that $M/M'$ is
a torsion $\phi$-module.  
\begin{conjecture}[The full Denis-Mordell-Lang conjecture]
\label{con:full-ml-dr}
Let $\phi_1:A\to K\{\tau\},\dots,\phi_g:A\to K\{\tau\}$ be Drinfeld
modules. Let $V\subset\mathbb{G}_a^g$ be an affine variety defined
over $\bK$. Let $\Gamma$ be a finite rank
$(\phi_1,\dots,\phi_g)$-submodule of $\mathbb{G}_a^g(\bK)$. Then there
exist algebraic $(\phi_1,\dots,\phi_g)$-submodules $B_1,\dots,B_l$ of
$\mathbb{G}_a^g$ and there exist $\gamma_1,\dots,\gamma_l\in\Gamma$
such that
$$V(\bK)\cap\Gamma = \bigcup_{i=1}^l (\gamma_i+B_i(\bK))\cap\Gamma.$$
\end{conjecture}
In \cite{Denis-conjectures}, Denis showed that under certain natural Galois theoretical assumptions, Conjecture~\ref{con:full-ml-dr} would follow from the weaker conjecture which would describe the intersection of an affine variety with a \emph{finitely generated $\phi$-module}.

Since then, the case $\Gamma$ is the product of the torsion submodules
of each $\phi_i$ was proved by Scanlon in \cite{Scanlon}, while
various other instances of Conjecture~\ref{con:full-ml-dr} were worked
out in \cite{IMRN} and \cite{full-ml-drinfeld}. We note that Denis
asked his conjecture also for \emph{$t$-modules}, which includes the
case of products of distinct Drinfeld modules acting on $\bG_a^g$.

For the sake of simplifying the notation, we denote by $\phi$ the action of $(\phi_1,\dots,\phi_g)$ on $\bG_a^g$.
We also note that if $V$ is an irreducible affine subvariety of
$\mathbb{G}_a^g$ which has a Zariski dense intersection with a finite
rank $\phi$-submodule $\Gamma$ of $\bG_a^g$, then the Denis-Mordell-Lang
conjecture predicts that $V$ is a translate of an algebraic
$\phi$-submodule of $\bG_a^g$ by a point in $\Gamma$. In particular,
if $V$ is an irreducible affine curve, which is \emph{not} a translate
of an algebraic $\phi$-submodule, then its intersection with any
finite rank $\phi$-submodule of $\bG_a^g$ should be finite.

In \cite{IMRN}, the first author studied the Denis-Mordell-Lang
conjecture for Drinfeld modules whose field of definition (for their
coefficients) is of transcendence degree at least equal to $2$. The
methods employed in \cite{IMRN} involve specializations, and so, it
was crucial for the $\phi$ there \emph{not} to be isomorphic with a
Drinfeld module defined over $\overline{\Fq(t)}$. In the present paper
we will study precisely this case left out in \cite{IMRN} and
\cite{full-ml-drinfeld}.  Our methods depend crucially on the
hypothesis that the transcendence degree of the field generated by the
coefficients of $\phi_i$ is one, since we use the fact that at each
place $v$, the number of residue classes in the ring of integers at
$v$ is finite.

The main result of our paper is describing the intersection of an
affine subvariety $V\subset\bG_a^g$ with a cyclic $\phi$-submodule
$\Gamma$ of $\bG_a^g$.

\begin{theorem}
\label{hypersurface}
Let $K$ be a function field of transcendence degree equal to one. Let $\phi_1:A\rightarrow K\{\tau\},\dots,\phi_g:A\to K\{\tau\}$ be
Drinfeld modules. Let
$(x_1,\dots,x_g)\in\bG_a^g(K)$ and let
$\Gamma\subset\mathbb{G}_a^g(K)$ be the cyclic
$(\phi_1,\dots,\phi_g)$-submodule generated by $(x_1,\dots,x_g)$. Let
$V\subset\mathbb{G}_a^g$ be an affine subvariety defined over $K$.
Then $V(K)\cap\Gamma$ is a finite union of cosets of
$(\phi_1,\dots,\phi_g)$-submodules of $\Gamma$.
\end{theorem}

Using an idea from \cite{full-ml-drinfeld}, we are able to extend the
above result to $(\phi_1,\dots,\phi_g)$-submodules of rank $1$ (see
our Theorem~\ref{full-rational}) in the special case where $V$ is a
curve.

\section{Proofs of our main results}
\label{proofs}

We continue with the notation from Section~\ref{notation}. Hence
$\phi_1,\dots,\phi_g$ are Drinfeld modules. We denote by $\phi$ the
action of $(\phi_1,\dots,\phi_g)$ on $\bG_a^g$. Also, let
$(x_1,\dots,x_g)\in\bG_a^g(K)$ and let $\Gamma$ be the cyclic
$\phi$-submodule of $\bG_a^g(K)$ generated by $(x_1,\dots,x_g)$.
Unless otherwise stated, $V\subset\bG_a^g$ is an affine subvariety
defined over $K$.

We first prove an easy combinatorial result which we will use in the proof of Theorem~\ref{hypersurface}.
\begin{lemma}
\label{cosets}
Let $\Gamma$ be a cyclic $\phi$-submodule of $\bG_a^g(K)$. Let
$\Gamma_0$ be a nontrivial $\phi$-submodule of $\Gamma$, and let
$S\subset\Gamma$ be an infinite set. Suppose that for every infinite
subset $S_0\subset S$, there exists a coset $C_0$ of $\Gamma_0$ such
that $C_0\cap S_0\ne\emptyset$ and $C_0\subset S$. Then $S$ is a
finite union of cosets of $\phi$-submodules of $\Gamma$.
\end{lemma}

\begin{proof}
  Since $S$ is infinite, $\Gamma$ is infinite, and thus $\Gamma$ is
  torsion-free.  Therefore, $\Gamma$ is an infinite cyclic
  $\phi$-module, which is isomorphic to $A$ (as a module over itself).
  Hence, via this isomorphism, $\Gamma_0$ is isomorphic to a
  nontrivial ideal $I$ of $A$. Since $A/I$ is finite (recall that
  $A=\Fq[t]$), there are finitely many cosets of $\Gamma_0$ in
  $\Gamma$. Thus, $S$ contains at most finitely many cosets of
  $\Gamma_0$.
  
  Now, let $\{y_i+\Gamma_0\}_{i=1}^\ell$ be all of the cosets of
  $\Gamma_0$ that are contained in $S$.  Suppose that
\begin{equation}
\label{wrong assumption}
S_0:=S\setminus\bigcup_{i=1}^\ell(y_i+\Gamma_0)\text{ is infinite.}
\end{equation}
Then using the hypotheses of this Lemma for $S_0$, we see that there
is a coset of $\Gamma_0$ that is contained in $S$ but is not one of
the cosets $(y_i+\Gamma_0)$ (because it has a nonempty intersection
with $S_0$). This contradicts the fact that
$\{y_i+\Gamma_0\}_{i=1}^\ell$ are {\it all} the cosets of $\Gamma_0$
that are contained in $S$.  Therefore $S_0$ must be finite.  Since any
finite subset of $\Gamma$ is a finite union of cosets of the trivial
submodule of $\Gamma$, this completes our proof.
\end{proof}

We will also use the following Lemma in the proof of Theorem~\ref{hypersurface}.
\begin{lemma}
\label{same ratio}
Let $\theta:A\to K\{\tau\}$ and $\psi:A\to K\{\tau\}$ be Drinfeld
modules. Let $v$ be a place of good reduction for $\theta$ and $\psi$.
Let $x,y\in\bC_v$. Let $0<r_v<1$, and let $B_v:=\{z\in\bC_v\mid
|z|_v<r_v\}$ be a sufficiently small ball centered at the origin with
the property that both $\log_{\theta,v}$ and $\log_{\psi,v}$ are
analytic isometries on $B_v$. Then for every polynomials $P, Q\in A$
such that $(\theta_P(x),\psi_P(y))\in B_v\times B_v$ and
$(\theta_Q(x),\psi_Q(y))\in B_v\times B_v$, we have
$$\log_{\theta,v}(\theta_P(x))\cdot\log_{\psi,v}(\psi_Q(y))=\log_{\theta,v}(\theta_Q(x))\cdot\log_{\psi,v}(\psi_P(y)).$$
\end{lemma}

\begin{proof}
  Since $v$ is a place of good reduction for $\theta$, all the
  coefficients of $\theta_Q$ are $v$-adic integers and thus,
  $|\theta_Q(\theta_P(x))|_v\le |\theta_P(x)|_v<r_v$ (we use the fact
  that $|\theta_P(x)|_v<r_v<1$, and so, each term of
  $\theta_Q(\theta_P(x))$ has its absolute value at most equal to
  $|\theta_P(x)|_v$). Using \eqref{logs and expos}, we conclude that
$$Q\cdot\log_{\theta,v}(\theta_P(x))=\log_{\theta,v}(\theta_{QP}(x))=\log_{\theta,v}(\theta_{PQ}(x))=P\cdot\log_{\theta,v}(\theta_Q(x)).$$
Similarly we obtain that $Q\cdot\log_{\psi,v}(\psi_P(x))=P\cdot\log_{\psi,v}(\psi_Q(x))$. This concludes the proof of Lemma~\ref{same ratio}.
\end{proof}

The following result is an immediate corollary of Lemma~\ref{same ratio}.
\begin{cor}
\label{same ratio 2}
With the notation as in Theorem~\ref{hypersurface}, assume in addition
that $x_1\notin\left(\phi_1\right)_{\tor}$. Let $v$ be a place of good
reduction for each $\phi_i$. Suppose $B_v$ is a small ball (of radius
less than $1$) centered at the origin such that each $\log_{\phi_i,v}$
is an analytic isometry on $B_v$. Then for each $i\in\{2,\dots,g\}$,
the fractions
$$\lambda_i:=\frac{\log_{\phi_i,v}\left(\left(\phi_i\right)_{P}(x_i)\right)}{\log_{\phi_1,v}\left(\left(\phi_1\right)_{P}(x_1)\right)}$$
are independent of the choice of the nonzero polynomial $P\in A$ for
which $\phi_P(x_1,\dots,x_g)\in B_v^g$.
\end{cor}

The following simple result on zeros of analytic functions can be
found in \cite[Proposition 2.1, p.  42]{Goss}.  We include a short
proof for the sake of completeness.

\begin{lemma}\label{zeros}
  Let $F(z) = \sum_{i=0}^{\infty} a_i z^i$ be a power series with
  coefficients in $\bC_v$ that is convergent in an open disc $B$ of
  positive radius around the point $z = 0$.  Suppose that $F$ is not
  the zero function.  Then the zeros of $F$ in $B$ are isolated.
\end{lemma}
\begin{proof}
  Let $w$ be a zero of $F$ in $B$.  We may rewrite $F$ in terms of
  $(z-w)$ as a power series $F(z) = \sum_{i=1}^{\infty} b_i (z-w)^i$
  that converges in a disc $B_w$ of positive radius around $w$. 
  Let $m$ be the smallest index $n$ such that $b_n \not= 0$.  

Because $F$ is convergent in $B_w$, then there exists a positive real number $r$ such that for all $n>m$, we have $\left|\frac{b_n}{b_m}\right|_v<r^{n-m}$. 
Then,
  for any $u \in B_w$ such that $0<|u-w|_v < \frac{1}{r}$, we
  have $|b_m (u-w)^m|_v > |b_n (u-w)^n|_v$ for all $n > m$. Hence $|F(u)|_v = |b_m(u-w)^m|_v\ne 0$. Thus $F(u)\ne 0$, and so, $F$ has no zeros other than $w$ in a nonempty
  open disc around $w$.
\end{proof}

We are ready to prove Theorem~\ref{hypersurface}.
\begin{proof}[Proof of Theorem~\ref{hypersurface}.]
We may assume $V(K)\cap\Gamma$ is infinite (otherwise the conclusion of Theorem~\ref{hypersurface} is obvisouly satisfied). Assuming $V(K)\cap\Gamma$ is infinite, we will show that there exists a nontrivial $\phi$-submodule $\Gamma_0\subset\Gamma$ such that each infinite subset of points $S_0$ in $V(K)\cap\Gamma$ has a nonempty intersection with a coset $C_0$ of $\Gamma_0$, and moreover, $C_0\subset V(K)\cap\Gamma$.
Then Lemma~\ref{cosets} will finish the proof of Theorem~\ref{hypersurface}.

First we observe that $\Gamma$ is not a torsion $\phi$-submodule. Otherwise $\Gamma$ is finite, contradicting our assumption that $V(K)\cap\Gamma$ is infinite. Hence, from now on, we assume (without loss of generality) that $x_1$ is not a torsion point for $\phi_1$.

We fix a finite set of polynomials $\{f_j\}_{j=1}^\ell \subset
K[X_1,\dots,X_g]$ which generate the vanishing ideal of $V$.

Let $v\in M_K$ be a place of $K$ which is of good reduction for all $\phi_i$ (for $1\le i\le g$). In addition, we assume each $x_i$ is integral at $v$ (for $1\le i\le g$). Then for each $P\in A$, we have 
$$\phi_P(x_1,\dots,x_g)\in \bG_a^g(\fO_v),$$
where $\fO_v$ is the ring
of $v$-adic integers in $K_v$ (the completion of $K$ at $v$). Because
$\fO_v$ is a compact space (we use the fact that $K$ is a function
field of transcendence degree $1$ and thus has a finite residue field
at $v$), we conclude that every infinite sequence of points
$\phi_P(x_1,\dots,x_g)\in V(K)\cap\Gamma$ contains a convergent
subsequence in $\fO_v^g$. Using Lemma~\ref{cosets}, it suffices to
show that there exists a nontrivial $\phi$-submodule $\Gamma_0\subset
\Gamma$ such that every convergent sequence of points in
$V(K)\cap\Gamma$ has a nonempty intersection with a coset $C_0$ of
$\Gamma_0$, and moreover, $C_0\subset V(K)\cap\Gamma$.

Now, let $S_0$ be an infinite subsequence of distinct points in
$V(K)\cap\Gamma$ which converges $v$-adically to
$(x_{0,1},\dots,x_{0,g})\in \fO_v^g$, let $0<r_v<1$, and let
$B_v:=\{z\in\bC_v\mid |z|_v<r_v\}$ be a small ball centered at the
origin on which each of the logarithmic functions $\log_{\phi_i,v}$ is
an analytic isometry (for $1\le i\le g$). Since
$(x_{0,1},\dots,x_{0,g})$ is the limit point for $S_0$, there exists a
$d\in A$ and an infinite subsequence $\{\phi_{d+P_n}\}_{n\ge 0}\subset
S_0$ (with $P_n=0$ if and only if $n=0$), such that for each $n\ge 0$,
we have
\begin{equation}
\label{hypersurface accumulation point}
\left|\left(\phi_i\right)_{d+P_n}(x_i) - x_{0,i}\right|_v <\frac{r_v}{2}
\quad \text{ for each $1\le i\le g$.}
\end{equation}
We will show that there exists an algebraic group $Y_0$, independent of $S_0$ and invariant under $\phi$, such that $\phi_d(x_1,\dots,x_g)+Y_0$ is a
subvariety of $V$ containing $\phi_{d+P_n}(x_1,\dots,x_g)$ for all
$P_n$. Thus the submodule $\Gamma_0:=Y_0(K)\cap\Gamma$ will satisfy the
hypothesis of Lemma~\ref{cosets} for the infinite subset
$V(K)\cap\Gamma\subset\Gamma$; this will yield the conclusion of
Theorem~\ref{hypersurface}.

Using \eqref{hypersurface accumulation point} for $n=0$ (we recall
that $P_0=0$), and then for arbitrary $n$, we see that
\begin{equation}
\label{hypersurface accumulation point near 0}
\left|\left(\phi_i\right)_{P_n}(x_i)\right|_v < \frac{r_v}{2}\text{ for each $1\le i\le g$.}
\end{equation}
Hence $\log_{\phi_i,v}$ is well-defined at
$\left(\phi_i\right)_{P_n}(x_i)$ for each $i\in\{1,\dots,g\}$ and for
each $n\ge 1$. Moreover, the fact that
$\left(\left(\phi_i\right)_{P_n+d}(x_i)\right)_{n\ge 1}$ converges to
a point in $\fO_v$ means that
$\left(\left(\phi_i\right)_{P_n}(x_i)\right)_{n\ge 1}$ converges to a
point which is contained in $B_v$ (see \eqref{hypersurface
  accumulation point near 0}).

Without loss of generality, we may assume 
\begin{equation}
\label{important}
|\log_{\phi_1,v}\left(\left(\phi_1\right)_{P_1}(x_1)\right)|_v = \max_{i=1}^g |\log_{\phi_i,v}\left(\left(\phi_i\right)_{P_1}(x_i)\right)|_v.
\end{equation}

Using the result of Corollary~\ref{same ratio 2}, we conclude that for
each $i\in\{2,\dots,g\}$, the following fraction is independent of $n$
and of the sequence $\{P_n\}_n$:
\begin{equation}
\label{hypersurface lambdas}
\lambda_i := \frac{\log_{\phi_i,v} \left(
    \left(\phi_i\right)_{P_n}(x_i)\right)}{\log_{\phi_1,v} \left(
    \left( \phi_1\right)_{P_n}(x_1)\right)}.
\end{equation}
Note that since $x_1$ is not a torsion point for $\phi_1$, the
denominator of $\lambda_i$ \eqref{hypersurface lambdas} is nonzero.
Because of equation \eqref{important}, we may conclude that
$|\lambda_i|_v\le 1$ for each $i$.

The fact that $\lambda_i$ is independent of the sequence $\{P_n\}_n$
will be used later to show that the $\phi$-submodule $\Gamma_0$ that
we construct is independent of the sequence $\{P_n\}_n$.

For each $n\ge 1$ and each $2\le i\le g$, we have
\begin{equation}
\label{hypersurface lambdas 2}
\log_{\phi_i,v} \left(\left(\phi_i\right)_{P_n}(x_i)\right) = \lambda_i\cdot\log_{\phi_1,v}\left(\left(\phi_1\right)_{P_n}(x_1)\right). 
\end{equation}
For each $i$, applying the exponential function $\exp_{\phi_i,v}$ to both sides of \eqref{hypersurface lambdas 2} yields
\begin{equation}
\label{hypersurface lambdas 3}
\left(\phi_i\right)_{P_n}(x_i) = \exp_{\phi_i,v}\left(\lambda_i\cdot\log_{\phi_1,v}\left(\left(\phi_1\right)_{P_n}(x_1)\right)\right). 
\end{equation}

Since $\phi_{d+P_n}\left(x_1,\dots,x_g\right)\in V(K)$, for each
$j\in\{1,\dots,\ell\}$ we have
\begin{equation}
\label{points on hypersurface}
f_j\left(\phi_{d+P_n}(x_1,\dots,x_{g})\right)=0\text{ for each $n$.}
\end{equation}

For each $j\in\{1,\dots,\ell\}$ we let $f_{d,j}\in K[X_1,\dots,X_{g}]$ be defined by 
\begin{equation}
\label{new hypersurface}
f_{d,j}\left(X_1,\dots,X_{g}\right):=f_j\left(\phi_d(x_1,\dots,x_{g})+(X_1,\dots,X_{g})\right).
\end{equation}
We let $V_d\subset\bG_a^g$ be the affine subvariety defined by the equations 
$$f_{d,j}(X_1,\dots,X_{g})=0\text{ for each $j\in\{1,\dots,\ell\}$.}$$
Using \eqref{points on hypersurface} and \eqref{new hypersurface}, we
see that for each $j\in\{1,\dots,\ell\}$ we have
\begin{equation}
\label{points on new hypersurface}
f_{d,j}\left(\phi_{P_n}(x_1,\dots,x_{g})\right)=0
\end{equation}
for each $n$, and so,
\begin{equation}
\label{points on new hypersurface 2}
\phi_{P_n}(x_1,\dots,x_g)\in V_d(K).
\end{equation}

For each $j\in\{1,\dots,\ell\}$, we let $F_{d,j}(u)$ be the analytic function defined on $B_v$ by
$$F_{d,j}(u):=f_{d,j}\left(u,\exp_{\phi_2,v}\left(\lambda_2\log_{\phi_1,v}(u)\right),\dots,\exp_{\phi_{g},v}\left(\lambda_{g}\log_{\phi_1,v}(u)\right)\right).$$
We note, because of \eqref{important}, and the fact that
$\log_{\phi_1,v}$ is an analytic isometry on $B_v$, that for each
$u\in B_v$, we have
\begin{equation}
\label{important 2}
|\lambda_i\cdot\log_{\phi_1,v}(u)|_v=|\lambda_i|_v\cdot |\log_{\phi_1,v}(u)|_v\le |u|_v<r_v.
\end{equation}
Equation \eqref{important 2} shows that $\lambda_i\cdot\log_{\phi_1,v}(u)\in B_v$, and so, $\exp_{\phi_i,v}\left(\lambda_i\cdot\log_{\phi_1,v}(u)\right)$ is well-defined.

Using \eqref{hypersurface lambdas 3} and \eqref{points on new hypersurface} we obtain that for every $n\ge 1$, we have
\begin{equation}
\label{hypersurface many zeros}
F_{d,j}\left(\left(\phi_1\right)_{P_n}(x_1)\right) = 0.
\end{equation}
Thus $\left(\left(\phi_1\right)_{P_n}(x_1)\right)_{n\ge 1}$ is a
sequence of zeros for the analytic function $F_{d,j}$ which has an
accumulation point in $B_v$. Lemma~\ref{zeros} then implies that
$F_{d,j}=0$, and so, for each $j\in\{1,\dots,\ell\}$, we have
\begin{equation}
\label{hypersurface morphism}
f_{d,j}\left(u,\exp_{\phi_2,v}\left(\lambda_2\log_{\phi_1,v}(u)\right),\dots,\exp_{\phi_{g},v}\left(\lambda_{g}\log_{\phi_1,v}(u)\right)\right) = 0
\end{equation}
For each $u\in B_v$, we let
$$Z_u := \left(u,\exp_{\phi_2,v}\left(\lambda_2\log_{\phi_1,v}(u)\right),\dots,\exp_{\phi_{g},v}\left(\lambda_{g}\log_{\phi_1,v}(u)\right)\right)\in\bG_a^g(\bC_v) .$$
Then \eqref{hypersurface morphism} implies that 
\begin{equation}
\label{Z_u is in Y_0}
Z_u\in V_d\text{ for each $u\in B_v$.}
\end{equation} 
Let $Y_0$ be the Zariski closure of $\{Z_u\}_{u\in B_v}$. Then
$Y_0\subset V_d$. Note that $Y_0$ is independent of the sequence
$\{P_n\}_n$ (because the $\lambda_i$ are independent of the sequence
$\{P_n\}_n$, according to Corollary~\ref{same ratio 2}).

We claim that for each $u\in B_v$ and for each $P\in A$, we have 
\begin{equation}
\label{Y_0 is invariant}
\phi_P(Z_u) = Z_{\left(\phi_1\right)_P(u)}.
\end{equation}
Note that for each $u\in B_v$, then also $(\phi_1)_P(u)\in B_v$ for
each $P\in A$, because each coefficient of $\phi_1$ is a $v$-adic
integer.  To see that \eqref{Y_0 is invariant} holds, we use
\eqref{logs and expos}, which implies that for each $i\in\{2,\dots,g\}$
we have
\begin{equation*}
\begin{split}  
  \exp_{\phi_i,v}\left(\lambda_i \log_{\phi_1,v}
    \left(\left(\phi_1\right)_P(u)\right)\right)
  & = \exp_{\phi_i,v}\left(\lambda_i\cdot P\cdot \log_{\phi_1,v}(u)\right)\\
  & =\exp_{\phi_i,v}\left(P\cdot\lambda_i\log_{\phi_1,v}(u)\right)\\
  & = \left(\phi_i\right)_P \left(
    \exp_{\phi_i,v}\left(\lambda_i\log_{\phi_1,v}(u)\right)\right).
\end{split}
\end{equation*}
Hence, \eqref{Y_0 is invariant} holds, and so, $Y_0$ is invariant
under $\phi$.  Furthermore, since all of the $\exp_{\phi_i,v}$ and
$\log_{\phi_i,v}$ are additive functions, we have
$Z_{u_1+u_2}=Z_{u_1}+Z_{u_2}$ for every $u_1,u_2\in B_v$.  Hence $Y_0$
is an algebraic group, which is also a $\phi$-submodule of $\bG_a^g$.  Moreover, $Y_0$ is
defined independently of $\Gamma$.

Let $\Gamma_0:=Y_0(K)\cap\Gamma$. Because $Y_0$ is invariant under $\phi$, then $\Gamma_0$ is a submodule of $\Gamma$. Because $Y_0\subset V_{d}$, it follows
that the translate $\phi_d(x_1,\dots,x_g)+Y_0$ is a subvariety of $V$
which contains all $\{\phi_{d+P_n}(x_1,\dots,x_g)\}_n$. In particular,
the (infinite) translate $C_0$ of $\Gamma_0$ by
$\phi_d(x_1,\dots,x_g)$ is contained in $V(K)\cap\Gamma$. Hence, every
infinite sequence of points in $V(K)\cap\Gamma$ has a nontrivial
intersection with a coset $C_0$ of (the nontrivial $\phi$-submodule)
$\Gamma_0$, and moreover, $C_0\subset V(K)\cap\Gamma$.  Applying
Lemma~\ref{cosets} thus finishes the proof of
Theorem~\ref{hypersurface}.
\end{proof}

In the course of our proof of Theorem~\ref{hypersurface} we also
proved the following statement.
\begin{theorem}
\label{stronger result}
Let $\Gamma$ be an infinite cyclic $\phi$-submodule of $\bG_a^g$. Then
there exists an infinite $\phi$-submodule $\Gamma_0\subset\Gamma$ such
that for every affine subvariety $V\subset\bG_a^g$, if
$V(\Kbar)\cap\Gamma$ is infinite, then $V(\Kbar)\cap\Gamma$ contains a
coset of $\Gamma_0$.
\end{theorem}

\begin{proof}
Let $v$ be a place of good reduction for $\phi$; in addition, we assume the points in $\Gamma$ are $v$-adic integers. Suppose that $V(\Kbar)\cap\Gamma$ is infinite. As shown in the proof of Theorem~\ref{hypersurface}, there exists a positive dimensional algebraic group $Y_0$, invariant under $\phi$, and depending only on $\Gamma$ and $v$ (but not on $V$), such that a translate of $Y_0$ by a point in $\Gamma$ lies in $V$. Moreover, $\Gamma_0:=Y_0(\Kbar)\cap\Gamma$ is infinite. Hence $\Gamma_0$ satisfies the conclusion of Theorem~\ref{stronger result}.
\end{proof}

\section{Further extensions}
\label{further extensions}

We continue with the notation from Section~\ref{proofs}:
$\phi_1,\dots,\phi_g$ are Drinfeld modules. As usual, we denote by
$\phi$ the action of $(\phi_1,\dots,\phi_g)$ on $\bG_a^g$. First we
prove the following consequence of Theorem~\ref{hypersurface}.
\begin{theorem}
\label{rational 2}
Let $V\subset\bG_a^g$ be an affine subvariety defined over $K$. Let
$\Gamma\subset\bG_a^g(K)$ be a finitely generated $\phi$-submodule of
rank $1$. Then $V(K)\cap\Gamma$ is a finite union of cosets of
$\phi$-submodules of $\Gamma$. In particular, if $V$ is an irreducible
curve which is not a translate of an algebraic $\phi$-submodule, then
$V(K)\cap\Gamma$ is finite.
\end{theorem}

\begin{proof}
  Since $A=\Fq[t]$ is a principal ideal domain, $\Gamma$ is the direct
  sum of its finite torsion submodule $\Gamma_{\tor}$ and a free
  submodule $\Gamma_1$, which is cyclic because $\Gamma$ has rank $1$. Therefore
  $$\Gamma = \bigcup_{\gamma\in\Gamma_{\tor}} \gamma+\Gamma_1,$$
  and
  so,
  $$V(K)\cap\Gamma = \bigcup_{\gamma\in\Gamma_{\tor}}
  V(K)\cap\left(\gamma + \Gamma_1\right) =
  \bigcup_{\gamma\in\Gamma_{\tor}} \left(\gamma+
    \left(-\gamma+V(K)\right)\cap\Gamma_1\right).$$
  Using the fact
  $\Gamma_{\tor}$ is finite and applying Theorem~\ref{hypersurface} to
  each intersection $\left(-\gamma+V(K)\right)\cap\Gamma_1$ thus completes our
  proof.
\end{proof}

We use the ideas from \cite{full-ml-drinfeld} to describe the
intersection of a curve $C$ with a $\phi$-module of rank $1$. So, let
$(x_1,\dots,x_g)\in\bG_a^g(K)$, let $\Gamma$ be the cyclic
$\phi$-submodule of $\bG_a^g(K)$ generated by $(x_1,\dots,x_g)$, and
let $\overline{\Gamma}$ be the $\phi$-submodule of rank $1$, containing
all $(z_1,\dots,z_g)\in\bG_a^g(\Kbar)$ for which there exists a
nonzero polynomial $P$ such that
$$\phi_P(z_1,\dots,z_g) \in\Gamma.$$
Since all polynomials $\phi_P$ (for $P\in A$) are separable, we have
$\overline{\Gamma}\subset\bG_a^g(K^{\sep})$.

With the notation above, we prove the following result; this may be
viewed as a Drinfeld module analog of McQuillan's result on semiabelian
varieties (see \cite{McQ}), which had been conjectured by Lang.
\begin{theorem}
\label{full-rational}
Let $C\subset\bG_a^g$ be an affine curve defined over $K$. Then $C(\Kbar)\cap\overline{\Gamma}$ is a finite union of cosets of $\phi$-submodules of $\overline{\Gamma}$.
\end{theorem}

Before proceeding to the proof of Theorem~\ref{full-rational} we first
prove two facts which will be used later. The first fact is an
immediate consequence of Theorem $1$ of \cite{Scanlon} (the
Denis-Manin-Mumford conjecture for Drinfeld modules), which we state below.

\begin{theorem}[Scanlon]
\label{Scanlon-MM}
Let $V\subset \mathbb{G}_a^g$ be an affine variety defined over
$\Kbar$. Then there exist algebraic $\phi$-submodules
$B_1,\dots,B_\ell$ of $\mathbb{G}_a^g$ and elements
$\gamma_1,\dots,\gamma_\ell$ of $\phi_{\tor}$ such that
$$V(\Kbar)\cap\phi_{\tor} = \bigcup_{i=1}^\ell \left(\gamma_i+B_i(\Kbar)\right)\cap\phi_{\tor}.$$
\end{theorem}

Moreover, in Remark $19$ from \cite{Scanlon}, Scanlon notes that his proof of
the Denis-Manin-Mumford conjecture yields a uniform bound on the
degree of the Zariski closure of $V(\Kbar)\cap\phi_{\tor}$, depending
only on $\phi$, $g$, and the degree of $V$. In particular, one obtains
the following uniform statement for translates of curves.

\begin{fact}
\label{C:needed uniformity}
Let $C\subset\mathbb{G}_a^g$ be an irreducible curve which is not a translate of an algebraic $\phi$-module of $\mathbb{G}_a^g$. Then there exists a positive integer $N$ such that for every $y\in\mathbb{G}_a^g(\Kbar)$, the set $\left(y+C(\Kbar)\right)\cap\phi_{\tor}$ has at most $N$ elements.
\end{fact}

\begin{proof}
  The curve $C$ contains no translate of a positive dimensional
  algebraic $\phi$-submodule of $\mathbb{G}_a^g$, so for every
  $y\in\mathbb{G}_a^g(\Kbar)$, the algebraic $\phi$-modules $B_i$
  appearing in the intersection
  $\left(y+C(\Kbar)\right)\cap\phi_{\tor}$ are all trivial. In
  particular, the set $\left(y+C(\Kbar)\right)\cap\phi_{\tor}$ is
  finite. Thus, using the uniformity obtained by Scanlon for his
  Manin-Mumford theorem, we conclude that the cardinality of
  $\left(y+C(\Kbar)\right)\cap\phi_{\tor}$ is uniformly bounded above
  by some positive integer $N$.
\end{proof}

We will also use the following fact in the proof of our Theorem~\ref{full-rational}.
\begin{fact}
\label{F:torsion}
Let $\phi:A\rightarrow K\{\tau\}$ be a Drinfeld module. Then for every positive integer $D$, there exist finitely many torsion points $y$ of $\phi$ such that $[K(y):K]\le D$.
\end{fact}

\begin{proof}
  If $y\in\phi_{\tor}$, then the canonical height $\hhat(y)$ of $y$
  (as defined in \cite{Denis}) equals $0$. Also, as shown in
  \cite{Denis}, the difference between the canonical height and the
  usual Weil height is uniformly bounded on $\Kbar$. Then
  Fact~\ref{F:torsion} follows by noting that there are finitely many
  points of bounded Weil height and bounded degree over the field $K$
  (using Northcott's theorem applied to the global function field $K$).
\end{proof}

We are now ready to prove Theorem~\ref{full-rational}.

\begin{proof}[Proof of Theorem~\ref{full-rational}.]
  Arguing as in the proof of Theorem~\ref{hypersurface}, it suffices
  to show that if $C$ is an irreducible affine curve (embedded in
  $\bG_a^g$), then $C(\Kbar)\cap\overline{\Gamma}$ is infinite only if
  $C$ is a translate of an algebraic $\phi$-submodule (because any translate of an algebraic $\phi$-module intersects $\overline{\Gamma}$ in a coset of some $\phi$-submodule of $\overline{\Gamma}$). Therefore, from
  now on, we assume $C$ is irreducible, that
  $C(\Kbar)\cap\overline{\Gamma}$ is infinite, and that $C$ is not a
  translate of an algebraic $\phi$-submodule.  We will derive a
  contradiction.

Let $z\in C(\Kbar)\cap\overline{\Gamma}$. For each field automorphism
$\sigma:K^{\sep}\to K^{\sep}$ that restricts to the identity on $K$,
we have $z^{\sigma}\in C\left(K^{\sep}\right)$ (because $C$ is defined
over $K$). By the definition of $\overline{\Gamma}$, there exists a
nonzero polynomial $P\in A$ such that $\phi_P(z)\in \Gamma$.  Since
$\phi_P$ has coefficients in $K$, we obtain
$$\phi_P\left(z^{\sigma}\right) = \left(\phi_P(z)\right)^{\sigma}=\phi_P(z).$$
The last equality follows from the fact that
$\phi_P(z)\in\Gamma\subset\bG_a^g(K)$.  We conclude that
$\phi_P\left(z^{\sigma}-z\right) =0$, and, thus, we have
$$T_{z,\sigma}:=z^{\sigma}-z\in\phi_{\tor}.$$
Moreover,
$T_{z,\sigma}\in (-z+C(\Kbar))\cap\phi_{\tor}$ (because $z^{\sigma}\in
C$). Using Fact~\ref{C:needed uniformity} we conclude that for each
fixed $z\in C(\Kbar)\cap\overline{\Gamma}$, the set
$\{T_{z,\sigma}\}_{\sigma}$ has cardinality bounded above by some
number $N$ (independent of $z$). In particular, this implies that $z$
has finitely many Galois conjugates, so $[K(z):K]\le N$.  Similarly we
have $\left[K(z^{\sigma}):K\right]\le N$; thus, we may conclude that
\begin{equation}
\label{torsion of bounded degree}
\left[K\left(T_{z,\sigma}\right):K\right]\le \left[K(z,z^{\sigma}):K\right]\le N^2.
\end{equation}
As shown by Fact~\ref{F:torsion}, there exists a finite set of torsion
points $w$ for which $[K(w):K]\le N^2$. Hence, recalling that $N$ is
independent of $z$, we see that the set
\begin{equation}
\label{few torsion}
H:=\{T_{z,\sigma}\}_{\substack{z\in C(\Kbar)\cap\overline{\Gamma}\\\sigma:K^{\sep}\to K^{\sep}}}\text{ is finite.}
\end{equation}

Now, since $H$ is a finite set of torsion points, there must exist a
nonzero polynomial $Q\in A$ such that $\phi_Q(H) =\{0\}$.
Therefore, $\phi_Q(z^{\sigma}-z)=0$ for each $z\in
C(\Kbar)\cap\overline{\Gamma}$ and each automorphism $\sigma$. Hence
$\phi_Q(z)^{\sigma} = \phi_Q(z)$ for each $\sigma$.  Thus, we have
\begin{equation}
\label{small extensions 2}
\phi_Q(z)\in \bG_a^g(K)\text{ for every $z\in C(\Kbar)\cap\overline{\Gamma}$.}
\end{equation}

Let $\Gamma_1:=\overline{\Gamma}\cap\bG_a^g(K)$.  Since
$\overline{\Gamma}$ is a finite rank $\phi$-module and $\bG_a^g(K)$ is
a \emph{tame} module (i.e. every finite rank submodule is finitely
generated; see \cite{Poonen} for a proof of this result), it follows
that $\Gamma_1$ is finitely generated. Let $\Gamma_2$ be the finitely
generated $\phi$-submodule of $\overline{\Gamma}$ generated by all
points $z\in\overline{\Gamma}$ such that $\phi_Q(z)\in\Gamma_1$. More
precisely, if $w_1,\dots,w_\ell$ generate the $\phi$-submodule
$\Gamma_1$, then for each $i\in\{1,\dots,\ell\}$, we find all the
finitely many $z_i$ such that $\phi_Q(z_i)=w_i$. Then this finite set
of all $z_i$ generate the $\phi$-submodule $\Gamma_2$. Thus $\Gamma_2$
is a finitely generated $\phi$-submodule, and moreover, using equation
\eqref{small extensions 2}, we obtain $C(\Kbar)\cap\overline{\Gamma} =
C(\Kbar)\cap\Gamma_2$.  Since $\Gamma_2$ is a finitely generated
$\phi$-submodule of rank $1$ (because
$\Gamma_2\subset\overline{\Gamma}$ and $\overline{\Gamma}$ has rank
$1$), Theorem~\ref{rational 2} finishes the proof of
Theorem~\ref{full-rational}.
\end{proof}

\def\cprime{$'$} \def\cprime{$'$} \def\cprime{$'$} \def\cprime{$'$}
\providecommand{\bysame}{\leavevmode\hbox to3em{\hrulefill}\thinspace}
\providecommand{\MR}{\relax\ifhmode\unskip\space\fi MR }
\providecommand{\MRhref}[2]{%
  \href{http://www.ams.org/mathscinet-getitem?mr=#1}{#2}
}
\providecommand{\href}[2]{#2}

\end{document}